\documentstyle[11pt,twoside]{article}
\include{graphicx}
\title{Exponentially small expansions associated with a generalised Mathieu series}
\author{\sc R. B.\ Paris\footnote{E-mail address: \ \ {\tt r.paris@abertay.ac.uk}} \\
{\em Division of Computing and Mathematics}, \\
{\em University of Abertay Dundee, Dundee DD1 1HG, UK}
}
\begin{document}
\def\f#1#2{\mbox{${\textstyle \frac{#1}{#2}}$}}
\def\dfrac#1#2{\displaystyle{\frac{#1}{#2}}}
\def\boldal{\mbox{\boldmath $\alpha$}}
\newcommand{\bee}{\begin{equation}}
\newcommand{\ee}{\end{equation}}
\newcommand{\lam}{\lambda}
\newcommand{\ka}{\kappa}
\newcommand{\al}{\alpha}
\newcommand{\om}{\omega}
\newcommand{\Om}{\Omega}
\newcommand{\fr}{\frac{1}{2}}
\newcommand{\fs}{\f{1}{2}}
\newcommand{\g}{\Gamma}
\newcommand{\br}{\biggr}
\newcommand{\bl}{\biggl}
\newcommand{\ra}{\rightarrow}
\newcommand{\gtwid}{\raisebox{-.8ex}{\mbox{$\stackrel{\textstyle >}{\sim}$}}}
\newcommand{\ltwid}{\raisebox{-.8ex}{\mbox{$\stackrel{\textstyle <}{\sim}$}}}
\renewcommand{\topfraction}{0.9}
\renewcommand{\bottomfraction}{0.9}
\renewcommand{\textfraction}{0.05}
\newcommand{\mcol}{\multicolumn}
\date{}
\maketitle
\begin{abstract}
We consider the generalised Mathieu series
\[\sum_{n=1}^\infty \frac{n^\gamma}{(n^\lambda+a^\lambda)^\mu}\qquad (\mu>0)\]
when the parameters $\lambda$ ($>0$) and $\gamma$ are even integers for large complex $a$ in the sector $|\arg\,a|<\pi/\lambda$. The
asymptotics in this case consist of a {\it finite} algebraic expansion together with an infinite sequence
of increasingly subdominant exponentially small expansions. When $\mu$ is also a positive integer it is possible to give closed-form evaluations of this series.
Numerical results are given to illustrate the accuracy of the expansion obtained.
\vspace{0.4cm}

\noindent {\bf Mathematics Subject Classification:} 30E15, 30E20, 34E05
\vspace{0.1cm}

\noindent {\bf Keywords:} asymptotic expansions, exponentially small expansions, generalised Mathieu series, Mellin transform method
\end{abstract}

\vspace{0.3cm}
\begin{center}
{\bf 1.\ Introduction}
\end{center}
\setcounter{section}{1}
\setcounter{equation}{0}
\renewcommand{\theequation}{\arabic{section}.\arabic{equation}}
This paper is a sequel to the asymptotic study of a generalised Mathieu series carried out by the author in \cite{P}. The functional series
\bee\label{e11}
\sum_{n=1}^\infty\frac{n}{(n^2+a^2)^\mu}
\ee
in the case $\mu=2$ 
was introduced by Mathieu in his 1890 book \cite{M} dealing with the elasticity of solid bodies. Considerable
effort has been devoted to the determination of upper and low bounds for the series with $\mu=2$ when the parameter $a>0$. Several integral representations for (\ref{e11}), and its alternating variant, have been obtained; see \cite{PST} and the references therein.

The asymptotic expansion of the more general functional series
\bee\label{e12}
S_{\mu,\gamma}(a;\lambda):=\sum_{n=1}^\infty\frac{n^\gamma}{(n^\lambda+a^\lambda)^\mu}\qquad (\mu>0,\ \lambda>0,\ \mu\lambda-\gamma>1)
\ee
was considered by Zvastavnyi \cite{Z} for $a\ra+\infty$ and by Paris \cite{P} for $|a|\ra\infty$ in the sector $|\arg\,a|<\pi/\lambda$. In \cite{P}, the additional factor $e_n:=\exp [-a^\lambda b/(n^\lambda+a^\lambda)]$ (with $b>0$) was included in the summand which, although not affecting the rate of convergence of the series (since $e_n\ra1$ as $n\ra\infty$), can modify the large-$a$ growth, particularly with the alternating variant of (\ref{e12}). Both these authors adopted a Mellin transform method and obtained the result\footnote{The restriction $\gamma>-1$ was imposed in \cite{P} to avoid the formation of a double pole for odd negative integer values of $\gamma$; in \cite{Z}, the parameter $\gamma$ was allowed to assume arbitrary real values.}, when $\gamma>-1$,
\bee\label{e120}
S_{\mu,\gamma}(a;\lambda)-\frac{\g(\frac{\gamma+1}{\lambda}) \g(\mu-\frac{\gamma+1}{\lambda})}{\lambda \g(\mu) a^{\lambda\mu-\gamma-1}}\sim\frac{1}{\g(\mu)}\sum_{k=0}^\infty \frac{(-)^k \g(\mu+k)}{k! a^{\lambda(k+\mu)}}\,\zeta(-\gamma-\lambda k)
\ee
as $|a|\ra\infty$ in the sector $|\arg\,a|<\pi/\lambda$, where $\zeta(s)$ denotes the Riemann zeta function.

In this paper we also employ the Mellin transform approach used in \cite{P, Z}, 
where our interest will be concerned with the parameter values $\mu>0$ and even integer values of $\lambda$ ($>0$) and $\gamma$. In this case, the asymptotic series on the right-hand side of (\ref{e120}) (the algebraic expansion) is either a finite series or vacuous on account of the trivial zeros of $\zeta(s)$. We shall find that the asymptotic expansion of $S_{\mu,\gamma}(a;\lambda)$ for large complex $a$ in the sector $|\arg\,a|<\pi/\lambda$ for these parameter values consists of a finite algebraic expansion together with an infinite sequence (when $\mu$ is not an integer) of increasingly subdominant exponentially small contributions. In the case of positive integer $\mu$ it is possible to give a closed-form evaluation of $S_{\mu,\gamma}(a;\lambda)$.

It is perhaps rather surprising that such an 
innocent-looking series should possess such an intricate
asymptotic structure. A similar phenomenon has been recently observed in the expansion of the generalised Euler-Jacobi series $\sum_{n=1}^\infty n^{-w}\exp\,[-an^p]$ as the parameter $a\ra 0$ when $p$ and $w$ are even integers; see \cite{P2}. 
The leading terms in the expansion of $S_{\mu,\gamma}(a;\lambda)$ when $\gamma=0$ and $\lambda=2, 4$ have been given in \cite{T} using the Poisson-Jacobi formula.

In the application of the Mellin transform method to the series in (\ref{e12}) and its alternating variant we shall require the following estimates for the gamma function and the Riemann zeta function. For real $\sigma$ and $t$, we have the estimates
\bee\label{e13}
\g(\sigma\pm it)=O(t^{\sigma-\frac{1}{2}}),\qquad |\zeta(\sigma\pm it)|=O(t^{\Omega(\sigma)} \log^\alpha t)\quad (t\ra+\infty),
\ee
where $\Omega(\sigma)=0$ ($\sigma>1$), $\fs-\fs\sigma$ ($0\leq\sigma\leq 1$), $\fs-\sigma$ ($\sigma<0$) and $\alpha=1$ ($0\leq\sigma\leq 1$), $\alpha=0$ otherwise \cite[p.~95]{ECT}. The zeta function $\zeta(s)$ has a simple pole of unit residue at $s=1$ and the evaluations for positive integer $k$
\[\zeta(0)=-\fs,\quad \zeta(-2k)=0,\quad \zeta(2k)=\frac{(2\pi)^{2k}}{2(2k)!}|B_{2k}|\ \ (k\geq 1),\]
\bee\label{e14}
B_0=1, \quad B_2=\f{1}{6}, \quad B_4=-\f{1}{30},\quad B_6=\f{1}{42}, \ldots ,
\ee
where $B_k$ are the Bernoulli numbers. Finally, we have the well-known functional relation satisfied by $\zeta(s)$ given by \cite[p.~603]{DLMF}
\bee\label{e15}
\zeta(s)=2^s \pi^{s-1} \zeta(1-s) \g(1-s) \sin \fs\pi s.
\ee
\vspace{0.6cm}

\begin{center}
{\bf 2.\ An integral representation}
\end{center}
\setcounter{section}{2}
\setcounter{equation}{0}
\renewcommand{\theequation}{\arabic{section}.\arabic{equation}}
The generalised Mathieu series defined in (\ref{e12}) can be written as
\bee\label{e21}
S_{\mu,\gamma}(a;\lambda)=a^{-\delta} \sum_{n=1}^\infty h(n/a),\qquad h(x):=\frac{x^\gamma}{(1+x^\lambda)^\mu},\quad \delta:=\lambda\mu-\gamma
\ee
where the parameter $\delta>1$
for convergence. We employ a Mellin transform approach as discussed in \cite[Section 4.1.1]{PK}. The Mellin transform of $h(x)$ is ${\cal H}(s)=\int_0^\infty x^{s-1}h(x)\,dx$, where
\begin{eqnarray*}
{\cal H}(s)&=&\int_0^\infty \frac{x^{\gamma+s-1}}{(1+x^\lambda)^\mu}\,dx=\frac{1}{\lambda}\int_0^\infty \frac{\tau^{(\gamma+s)/\lambda-1}}{(1+\tau)^\mu}\,d\tau\\
\\
&=&\frac{\g(\frac{\gamma+s}{\lambda}) \g(\mu-\frac{\gamma+s}{\lambda})}{\lambda \g(\mu)}
\end{eqnarray*}
in the strip $-\gamma<\Re (s)<\delta$. Using the Mellin inversion theorem (see, for example, \cite[p.~118]{PK}), we find
\[\sum_{n=1}^\infty h(n/a)=\frac{1}{2\pi i}\sum_{n=1}^\infty \int_{c-\infty i}^{c+\infty i} {\cal H}(s) (n/a)^{-s}ds=\frac{1}{2\pi i}\int_{c-\infty i}^{c+\infty i} {\cal H}(s) \zeta(s) a^s ds,\]
where $1<c<\delta$. The inversion of the order of summation and integration is justified by absolute convergence provided $1<c<\delta$. 

Then, from (\ref{e21}), we have \cite{P, Z}
\bee\label{e22}
S_{\mu,\gamma}(a;\lambda)=\frac{a^{-\delta}}{\lambda \g(\mu)}\,\frac{1}{2\pi i}\int_{c-\infty i}^{c+\infty i}
\g\bl(\frac{\gamma+s}{\lambda}\br) \g\bl(\mu-\frac{\gamma+s}{\lambda}\br) \zeta(s)\,a^s ds,
\ee
where $1<c<\delta$. From the estimates in (\ref{e13}), the integral in (\ref{e22}) then defines $S_{\mu,\gamma}(a;\lambda)$ for complex $a$ in the sector $|\arg\,a|<\pi/\lambda$. The asymptotic expansion of
$S_{\mu,\gamma}(a;\lambda)$ for large $a$ and real parameters $\lambda$, $\mu$ and $\gamma$, such that $\delta>1$ is given in (\ref{e13}); see \cite[Theorem 3]{P} for complex $a$ and \cite{Z} for positive $a$ and unrestricted $\gamma$. 
 
We now suppose in the remainder of this paper that $\mu>0$, with $\lambda$ ($>0$) and $\gamma$ both chosen to be even integers. More specifically, we write
\bee\label{e23}
\lambda=2p\ \ (p=1, 2, \ldots ),\quad \gamma=2m \ \ (m=0,\pm 1, \pm 2, \ldots ).
\ee
The integration path in (\ref{e22}) lies to the right of the simple pole of $\zeta(s)$ at $s=1$ and the poles of $\g((\gamma+s)/\lambda)$ at $s=-\gamma-\lambda k$ ($k=0, 1, 2, \ldots\,$), but to the left of the poles of the second gamma function at $s=\delta+\lambda k$. When $\gamma=2, 4, \ldots\,$, the poles at $s=-\gamma-\lambda k$ are cancelled by the trivial zeros of $\zeta(s)$ at $s=-2, -4, \ldots $. When $\gamma=-2m$, $m=0, 1, 2, \ldots\,$, however, there is
a finite set of poles of this sequence  on the left of the integration path situated in $\Re (s)\geq 0$ with $0\leq k\leq k^*$, where $k^*$ is the index that satisfies 
\bee\label{e23a}
m-k^*p\geq 0,\qquad m-(k^*+1)p<0.
\ee
The remaining poles of this sequence corresponding to $k>k^*$ are cancelled by the trivial zeros of $\zeta(s)$ with the result that there are again no poles in $\Re (s)<0$.

We consider the integral in (\ref{e22}) taken round the rectangular contour with vertices at $c\pm iT$ and $-c'\pm iT$,
where $c'>0$ and $T>0$. The contribution from the upper and lower sides of the rectangle $s=\sigma\pm iT$, $-c'\leq\sigma\leq c$, vanishes as $T\ra\infty$ provided $|\arg\,a|<\pi/\lambda$, since from (\ref{e13}), the modulus of the integrand is controlled by $O(T^{\Omega(\sigma)+\frac{1}{2}\mu-\frac{3}{4}} \log\,T\,e^{-\Delta T})$, where $\Delta=\pi/\lambda-|\arg\,a|$. Evaluation of the residues then yields
\bee\label{e24}
S_{\mu,\gamma}(a;\lambda)=\frac{\g(\frac{\gamma+1}{\lambda}) \g(\mu-\frac{\gamma+1}{\lambda})}{\lambda \g(\mu) a^{\delta-1}}+H_{\mu,\gamma}(a;\lambda)+J(a),
\ee
where the finite algebraic expansion $H_{\mu,\gamma}(a;\lambda)$ (with $\gamma=2m$) is given by
\bee\label{e25}
H_{\mu,\gamma}(a;\lambda)=\left\{\begin{array}{ll} \!\!0 & \ \ (m=1, 2, \ldots) \\
\displaystyle{
\!\!\frac{a^{-\lambda\mu}}{\g(\mu)}\sum_{k=0}^{k^*}\frac{(-)^k \g(\mu+k)}{k! \,a^{\lambda k}}\,\zeta(2m-2kp)} & \ \ (m=0, -1, -2, \ldots), \end{array}\right.
\ee
with the index $k^*$ being defined in (\ref{e23a}), and
\bee\label{e26}
J(a)=\frac{a^{-\delta}}{\lambda \g(\mu)}\,\frac{1}{2\pi i}\int_{-c'-\infty i}^{-c'+\infty i}
\g\bl(\frac{\gamma+s}{\lambda}\br) \g\bl(\mu-\frac{\gamma+s}{\lambda}\br) \zeta(s)\,a^s ds\quad (c'>0).
\ee
The values of $\zeta(s)$ at $s=0, 2, 4, \ldots$ can be expressed in terms of the Bernoulli numbers, if so desired, by (\ref{e14}).

The integrand in $J(a)$ is holomorphic in $\Re (s)<0$, so that further displacement of the contour to the left can produce no additional terms in the algebraic expansion
of $S_{\mu,\gamma}(a;\lambda)$. We shall see in the next section that $J(a)$ possesses an infinite sequence of increasingly exponentially small terms in the large-$a$ limit.
\vspace{0.6cm}

\begin{center}
{\bf 3.\ The exponentially small expansion of $J(a)$}
\end{center}
\setcounter{section}{3}
\setcounter{equation}{0}
\renewcommand{\theequation}{\arabic{section}.\arabic{equation}}
In the integral (\ref{e26}), we make the change of variable $s\ra -s-\gamma$ to find
\[J(a)=\frac{a^{-\lambda\mu}}{\lambda \g(\mu)}\,\frac{1}{2\pi i}\int_{d-\infty i}^{d+\infty i}
\g\bl(\frac{-s}{\lambda}\br) \g\bl(\mu+\frac{s}{\lambda}\br) \zeta(-s-\gamma)\,a^{-s} ds,\qquad d=c'-\gamma.\]
We now employ (\ref{e13}) to convert the zeta function into one with real part greater than unity. With the parameters $\lambda$ and $\gamma$ in (\ref{e23}), the above integral can then be written in the form
\[J(a)=\frac{(-)^m a^{-\lambda\mu}}{(2\pi)^{\gamma}\lambda \g(\mu)}\,\frac{1}{2\pi i}\int_{d-\infty i}^{d+\infty i}
G(s)\,\zeta(1\!+\!s\!+\!\gamma) (2\pi a)^{-s}\,\frac{\sin (\fs\pi s)}{\sin (\frac{\pi s}{\lambda})}\,ds,\]
where
\bee\label{e300}
G(s):=\frac{\g(s\!+\!\gamma\!+\!1) \g(\mu\!+\!\frac{s}{\lambda})}{\g(1\!+\!\frac{s}{\lambda})}.
\ee

Making use of the expansion
\bee\label{e30}
\frac{\sin (\fs\pi s)}{\sin (\frac{\pi s}{\lambda})}\equiv\frac{\sin (\frac{\pi ps}{\lambda})}{\sin (\frac{\pi s}{\lambda})}=\sum_{r=0}^{p-1} e^{-i\omega_rs},\qquad \omega_r:=(p\!-\!1\!-\!2r)\frac{\pi}{\lambda},
\ee
we then obtain
\bee\label{e31}
J(a)=\frac{(-)^m a^{-\lambda\mu}}{(2\pi)^{\gamma}\lambda \g(\mu)}\,\sum_{r=0}^{p-1} {\cal E}_r(a),
\ee
where 
\bee\label{e32}
{\cal E}_r(a)=\frac{1}{2\pi i}\int_{d-\infty i}^{d+\infty i}G(s)\,\zeta(s\!+\!\gamma\!+\!1)\,(2\pi ae^{ i\omega_r})^{-s}ds
\ee
and $d+\gamma=c'>0$.

The integrals ${\cal E}_r(a)$ ($0\leq r\leq p-1$) have no poles to the right of the integration path, so that we can displace the path as far to the right as we please. On such a displaced path, which we denote by  $L$, $|s|$ is everywhere large. Let $M$ denote an arbitrary positive integer. The quotient of gamma functions in $G(s)$ may be expanded by appealing to the inverse-factorial expansion given in \cite[p.~53]{PK} to obtain
\bee\label{e33}
G(s)=\lambda^{1-\mu} \bl\{\sum_{j=0}^{M-1}(-)^j c_j \g(s+\vartheta-j)+\rho_M(s) \g(s+\vartheta-M)\br\},\quad \vartheta:=\mu+\gamma,
\ee
where $c_0=1$ and $\rho_M(s)=O(1)$ as $|s|\ra\infty$ in $|\arg\,s|<\pi$. An algorithm for the evaluation of the coefficients $c_j$ is discussed in Section 4.
Substitution of the expansion (\ref{e33}) into (\ref{e32}) then produces
\bee\label{e34}
{\cal E}_r(a)=\lambda^{1-\mu}\bl\{\sum_{j=0}^{M-1}(-)^jc_j\,\frac{1}{2\pi i}\int_L\g(s\!+\!\vartheta\!-\!j) \zeta(s\!+\!\gamma\!+\!1)\,(2\pi ae^{i\omega_r})^{-s}ds+R_{M,r}\br\},
\ee
where the remainders $R_{M,r}$ are given by
\bee\label{e34a}
R_{M,r}=\frac{1}{2\pi i}\int_L\rho_M(s) \g(s\!+\!\vartheta\!-\!M) \zeta(s\!+\!\gamma\!+\!1)\,(2\pi ae^{i\omega_r})^{-s}ds.
\ee

The integrals appearing in (\ref{e34}) can be evaluated by making use of the well-known result 
\[\frac{1}{2\pi i}\int_{L'} \g(s+\alpha) z^{-s}ds=z^\alpha e^{-z}\qquad (|\arg\,z|<\fs\pi),\]
where $L'$ is a path parallel to the imaginary $s$-axis lying to the right of all the poles of $\g(s+\alpha)$;
see, for example, \cite[Section 3.3.1]{PK}. Upon expansion of the zeta function in (\ref{e34}) (since on $L$ its argument satisfies $\Re (s)+\gamma+1>1$) we find
\[\frac{1}{2\pi i}\int_L\g(s\!+\!\vartheta\!-\!j) \zeta(s\!+\!\gamma\!+\!1)\,(2\pi ae^{i\omega_r})^{-s}ds=\sum_{n=1}^\infty \frac{(2\pi nae^{i\omega_r})^{\vartheta-j}}{n^{1+\gamma}} \exp\,[-2\pi nae^{i\omega_r}] \]
\bee\label{e35a}
=X_r^{\vartheta-j} e^{-X_r} K_j(X_r;\mu),\qquad X_r:=X e^{i\omega_r},\quad X:=2\pi a,
\ee
where we have defined the exponential sum
\bee\label{e35}
K_j(X_r;\mu):=\sum_{n=1}^\infty \frac{e^{-(n-1)X_r}}{n^{1-\mu+j}}.
\ee
This evaluation is valid provided that the variable $X_r$ satisfies the convergence conditions
\[|\arg\,a+\omega_r|<\fs\pi\qquad (0\leq r\leq p-1).\]
From the definition of $\omega_r$ in (\ref{e30}), it is easily verified that these conditions are met when $|\arg\,a|<\pi/\lambda$. It is then evident that $K_j(X_r;\mu)\sim 1$ as $|a|\ra\infty$ in $|\arg\,a|<\pi/\lambda$.

Thus we find
\bee\label{e36}
{\cal E}_r(a)=\lambda^{1-\mu}e^{-X_r} \sum_{j=0}^{M-1}(-)^j c_j X_r^{\vartheta-j}\,K_j(X_r;\mu)+R_{M,r}.
\ee
Bounds for the remainders of the type $R_{M,r}$ have been considered in \cite[p.~71]{PK}; see also \cite[\S 10.1]{Br}. The integration path in (\ref{e34a}) is such that $\Re (s)+\gamma+1>1$, so that we may employ the bound $|\zeta(x+iy)|\leq\zeta(x)$ for real $x$, $y$ with $x>1$. A slight modification of Lemma 2.7 in \cite[p.~71]{PK} then shows that
\bee\label{e37}
R_{M,r}=O(X_r^{\vartheta-M} e^{-X_r})
\ee
as $|a|\ra\infty$ in the sector $|\arg\,a|<\pi/\lambda$.

The expansion of $S_{\mu,\gamma}(a;\lambda)$ then follows from (\ref{e24}), (\ref{e31}), (\ref{e36}) and (\ref{e37})
and is given in the following theorem.
\newtheorem{theorem}{Theorem}
\begin{theorem}$\!\!\!.$
Let $\mu>0$, $\gamma=2m$, $\lambda=2p$, where $m=0, \pm1, \pm2, \ldots$ and $p=1, 2, \ldots\,$. Further, let $M$ denote a positive integer, $\omega_r=\pi(p-1-2r)/(2p)$ for $0\leq r\leq p-1$ and $\delta=\lambda\mu-\gamma$, $\vartheta=\mu+\gamma$. Then 
\bee\label{e38}
S_{\mu,\gamma}(a;\lambda)=\frac{\g(\frac{\gamma+1}{\lambda}) \g(\mu-\frac{\gamma+1}{\lambda})}{\lambda \g(\mu) a^{\delta-1}}+H_{\mu,\gamma}(a;\lambda)+\frac{(-)^m}{\g(\mu)} \bl(\frac{\pi}{p}\br)^\mu a^{\mu-\delta} \sum_{r=0}^{p-1} E_r(a)
\ee
as $|a|\ra\infty$ in the sector $|\arg\,a|<\pi/\lambda$. The finite algebraic expansion $H_{\mu,\gamma}(a;\lambda)$ is defined in (\ref{e25}) and the exponentially small expansions $E_r(a)$ are given by
\bee\label{e38a}
E_r(a)=e^{-X_r+i\vartheta\omega_r}\bl\{\sum_{j=0}^{M-1}(-)^jc_j X_r^{-j} K_j(X_r;\mu)+O(X_r^{-M})\br\}\qquad (0\leq r\leq p-1),
\ee
where the leading coefficient $c_0=1$ and $X_r=2\pi ae^{i\omega_r}$. The infinite exponential sums $K_j(X_r;\mu)$ are defined in (\ref{e35}).
\end{theorem}

When $a$ is a real variable, the expansion in Theorem 1 can be expressed in a different form. We have the following theorem.
\begin{theorem}$\!\!\!.$
Let the parameters $\mu$, $\gamma$, $\lambda$ and $\delta$, $\vartheta$, $\omega_r$ be as in Theorem 1. Then with $N=\lfloor\fs p\rfloor$ and $X=2\pi a$, the expansion for $S_{\mu,\gamma}(a;\lambda)$ becomes
\[S_{\mu,\gamma}(a;\lambda)=\frac{\g(\frac{\gamma+1}{\lambda}) \g(\mu-\frac{\gamma+1}{\lambda})}{\lambda \g(\mu) a^{\delta-1}}+H_{\mu,\gamma}(a;\lambda)\hspace{6cm}\]
\bee\label{e39}
\hspace{3cm}+\frac{(-)^m}{\g(\mu)}\bl(\frac{\pi}{p}\br)^\mu a^{\mu-\delta}
\bl\{\sum_{r=0}^{N-1} E_r^*(a)+ \bl(\!\!\begin{array}{c}0\\ \fs E_N^*(a)\end{array}\!\!\br) \br\}\qquad \bl\{\!\!\begin{array}{l} p\ \mbox{even}\\p\ \mbox{odd}\end{array}
\ee
as $a\ra+\infty$, where for arbitrary positive integer $M$
\bee\label{e39a}
E_r^*(a)=2e^{-X\cos \omega_r} \bl\{\sum_{j=0}^{M-1}(-)^jc_jX^{-j} K_j^*(X;\omega_r)+O(X^{-M})\br\}\qquad (0\leq r\leq N)
\ee
and the infinite exponential sums $K^*_j(X;\omega_r)$ are defined by
\[K_j^*(X;\omega_r)=\sum_{n=1}^\infty\frac{e^{-(n-1)X\cos \omega_r}}{n^{1-\mu+j}}\,\cos\,[nX\sin \omega_r+(j-\vartheta)\,\omega_r].\]
When $p$ is odd, the quantity $\omega_N=0$.
\end{theorem}

\vspace{0.6cm}

\begin{center}
{\bf 4.\ The coefficients $c_j$}
\end{center}
\setcounter{section}{4}
\setcounter{equation}{0}
\renewcommand{\theequation}{\arabic{section}.\arabic{equation}}
We describe an algorithm for the computation of the coefficients $c_j$ that appear in the exponentially small expansions $E_r(a)$ and $E_r^*(a)$ in (\ref{e38a}) and (\ref{e39a}).
The expression for the ratio of gamma functions in $G(s)$ in (\ref{e33}) may be written in the form
\[\frac{G(s)}{\g(s+\vartheta)}=\lambda^{1-\mu}\bl\{\sum_{j=0}^{M-1} \frac{c_j}{(1- s-\vartheta)_j}+\frac{\rho_M(s)}{(1-s-\vartheta)_M}\br\},\]
where $(\alpha)_j=\g(\alpha+j)/\g(\alpha)$ is the Pochhammer symbol. If we introduce the scaled gamma function $\g^*(z)=\g(z)/(\sqrt{2\pi}\,z^{z-\fr}e^{-z})$, 
then we have
\[\g(\beta s+\gamma)=\g^*(\beta s+\gamma) (2\pi)^\fr e^{-\beta s} (\beta s)^{\beta s+\gamma-\fr}\,{\bf e}(\beta s;\gamma),\]
where
\[{\bf e}(\beta s;\gamma):= \exp \bl[(\beta s+\gamma-\fs) \log (1+\frac{\gamma}{\beta s})-\gamma\br].\]
The above ratio of gamma functions may therefore be expressed as
\bee\label{e41}
R(s)\Upsilon(s)=\sum_{j=0}^{M-1}\frac{c_j}{(1-s-\vartheta)_j}+\frac{\rho_M(s)}{(1-s-\vartheta)_M}
\ee
as $|s|\ra\infty$ in $|\arg\,s|<\pi$, where
\[R(s):=\frac{{\bf e}(s;\gamma\!+\!1)\,{\bf e}(s/\lambda; \mu)}{{\bf e}(s/\lambda; 1)\,{\bf e}(s;\vartheta)},\qquad \Upsilon(s):=\frac{\g^*(s\!+\!\gamma\!+\!1)\,\g^*(\mu\!+\!s/\lambda)}{\g^*(1\!+\!s/\lambda)\,\g^*(s\!+\!\vartheta)}~. \]

We now let $\xi:=s^{-1}$ and follow the procedure described in \cite[p.47]{PK}. We expand $R(s)$ and $\Upsilon(s)$ for $\xi\ra 0$ making use of the well-known expansion (see, for example, \cite[p.~71]{PK})
\[\g^*(z)\sim\sum_{k=0}^\infty(-)^k\gamma_kz^{-k}\qquad(|z|\ra\infty;\ |\arg\,z|<\pi),\]
where $\gamma_k$ are the Stirling coefficients, with 
\[\gamma_0=1,\quad \gamma_1=-\f{1}{12},\quad \gamma_2=\f{1}{288},\quad  \gamma_3=\f{139}{51840},
\quad \gamma_4=-\f{571}{2488320},\, \ldots\ .\]
After some straightforward algebra we find that
\[R(s)=1+\fs(\mu-1)\{(\lambda-1)\mu-2\gamma\}\xi+O(\xi^2), \]
\[\Upsilon(s)=1-\f{1}{12}(\mu-1)(\lambda^2-1)\xi^2+O(\xi^3),\]
so that upon equating coefficients of $\xi$ in (\ref{e41}) we can obtain $c_1$.
The higher coefficients can be obtained
by matching coefficients recursively with the aid of {\it Mathematica} to find 
\[c_0=1,\quad c_1=\f{1}{2}(\mu-1)\{2\gamma-(\lambda-1)\mu\},\]
\[c_2=\f{1}{24}(\mu-1)(\mu-2)\{12\gamma(\gamma-(\lambda-1)\mu-1)+(\lambda-1)\mu(5-3\mu+\lambda(3\mu-1))\},\]
\[c_3=-\f{1}{48}(\mu-1)(\mu-2)(\mu-3)\{2-2\gamma+(\lambda-1)\mu\}\{4\gamma(\gamma-(\lambda-1)\mu-2)\]
\bee\label{e42}
\hspace{5cm}+(\lambda-1)\mu(3+\lambda(\mu-1)-\mu)\},\ldots\ .
\ee
The rapidly increasing complexity of the coefficients with $j\geq 4$ prevents their presentation. However,
this procedure is found to work well in specific cases when the various parameters have numerical values, where many coefficients have been so calculated.
In Table 1 we present some values\footnote{In the tables we write the values as $x(y)$ instead of $x\times 10^y$.} of the coefficients $c_j$ for $1\leq j\leq 10$, which are used in the specific examples considered in Section 5.
\begin{table}[t]
\caption{\footnotesize{The coefficients $c_j$ ($1\leq j\leq 10$) for different $\gamma$ when $\mu=5/4$ and $\lambda=4$.}}
\begin{center}
\begin{tabular}{|r|c|c|c|}
\hline
&&&\\[-0.3cm]
\mcol{1}{|c|}{$j$} & \mcol{1}{c|}{$\gamma=0$} &\mcol{1}{c|}{$\gamma=2$}& \mcol{1}{c|}{$\gamma=-2$} \\
[.1cm]\hline
&&&\\[-0.25cm]
1 & $-4.6875000000(-1)$ & $+3.1250000000(-2)$ & $-2.1250000000(+0)$ \\
2 & $-3.5888671875(-1)$ & $+1.5673828125(-1)$ & $-2.5546875000(+0)$ \\
3 & $-4.0534973145(-1)$ & $+2.3551940918(-1)$ & $-6.8701171875(+0)$ \\
4 & $-3.3581793308(-1)$ & $+1.6646325588(-1)$ & $-2.5683746338(+1)$ \\
5 & $+7.5268601999(-1)$ & $-9.0884858742(-1)$ & $-1.1944799423(+2)$ \\
6 & $+6.4821335676(+0)$ & $-6.6501405553(+0)$ & $-6.6193037868(+2)$ \\
7 & $+2.6358910987(+1)$ & $-2.7627888119(+1)$ & $-4.2794038211(+3)$ \\
8 & $+4.5855530043(+1)$ & $-5.8401959193(+1)$ & $-3.1831413077(+4)$ \\
9 & $-3.7955573596(+2)$ & $+2.8858407940(+2)$ & $-2.6901844936(+5)$ \\
10& $-5.1286970180(+3)$ & $+4.6231064924(+3)$ & $-2.5504879368(+6)$ \\
[.2cm]\hline
\end{tabular}
\end{center}
\end{table}
\vspace{0.3cm}

\noindent{4.1}\ \ {\it The coefficients $c_j$ when $\lambda=2$, $\gamma=0$}
\vspace{0.3cm}

\noindent
When $\lambda=2$ and $\gamma=0$, it is possible to express the coefficients $c_j$ in closed form for arbitrary $\mu>0$.
From (\ref{e300}), we have
\[G(s)=\frac{\g(1+s) \g(\mu+\fs s)}{\g(1+\fs s)}=\frac{2^s}{\sqrt{\pi}}\,\g(\fs+\fs s) \g(\mu+\fs s)\]
upon use of the duplication formula for the gamma function. The inverse factorial expansion of a product of two gamma functions with equal coefficients of $s$ is given in \cite[pp.~51--52]{PK} in the form
\[\g(s+\alpha) \g(s+\beta)\sim 2^{\frac{3}{2}-\alpha-\beta-2s}\sqrt{\pi} \sum_{j=0}^\infty d_j\, \g(2s\!+\!\alpha\!+\!\beta\!-\!\fs\!-\!j)\]
as $|s|\ra\infty$ in $|\arg\,s|<\pi$,
where the coefficients satisfy $d_0=1$ and
\[d_j=\frac{2^{-j}}{j!} \prod_{r=1}^j \{(\alpha-\beta)^2-(r-\fs)^2\}\qquad (j\geq 1).\]

Putting $\alpha=\fs$ and $\beta=\mu$, with $s\ra \fs s$, we therefore obtain the coefficients in the inverse factorial expansion of $G(s)$ when $\lambda=2$, $\gamma=0$ and $\mu>0$ given by
\bee\label{e42a}
c_j=\frac{(-2)^{-j}}{j!}\,\prod_{r=1}^j (\mu-r)(\mu+r-1)\qquad (j\geq 1).
\ee

\vspace{0.3cm}

\noindent{4.2}\ \ {\it The coefficients $c_j$ when $\mu$ is an integer}
\vspace{0.3cm}

\noindent
A study of the coefficients $c_j$ with the aid of {\it Mathematica} enables us to conjecture that they possess the general form
\[c_j=(\mu-1)(\mu-2) \ldots (\mu-j)\,P_j(\mu,\gamma,\lambda)\qquad (j\geq 1),\]
where $P_j$ denotes a polynomial of degree $j$ in the parameters $\mu$, $\gamma$ and $\lambda$. This implies that the sequence of coefficients is finite for integer values of $\mu$; that is, for positive integer $q$, we have
\[c_j=0\qquad (j\geq q;\ \mu=q,\  q=1, 2, \ldots).\] 
This can also be seen from (\ref{e300}) where, with $\mu=q$,
\bee\label{e43a}
G(s)=\frac{\g(s\!+\!\gamma\!+\!1)\g(q\!+\!\frac{s}{\lambda}) }{\g(1\!+\!\frac{s}{\lambda})}
=\g(s\!+\!\gamma\!+\!1)\,\bl(1+\frac{s}{\lambda}\br)_{q-1}.
\ee
When $q=1$, we have $\vartheta=1+\gamma$ and the expansion (\ref{e33}) is satisfied trivially by terminating the series at the leading term with $\rho_1(s)\equiv 0$. When $q\geq 2$, the right-hand side of (\ref{e33}) must terminate at $M=q$ with $\rho_q(s)\equiv 0$ to yield
\bee\label{e43}
G(s)=\lambda^{1-q} \g(s\!+\!\gamma\!+\!1) \sum_{j=0}^{q-1}(-)^j c_j\,(s+\gamma+1)_{g-j-1}
\ee
in order to have the polynomials in $s$ in (\ref{e43a}) and (\ref{e43}) of the same degree.

Then the exponential expansions $E_r(a)$ in (\ref{e38a}) become the {\it finite} sums
\[E_r(a)=e^{-X_r+i\vartheta\omega_r} \sum_{j=0}^{q-1} (-)^jc_j X_r^{-j} K_j(X_r;q),\qquad\vartheta=q+\gamma,\]
where the infinite sums $K_j(X_r;q)$ may be expressed exactly in terms of derivatives of an exponential by
\begin{eqnarray}
K_j(X_r;q)&=&e^{X_r}\sum_{n=1}^\infty n^{q-j-1} e^{-nX_r}=(-1)^{q-j-1}e^{X_r}D^{q-j-1} \sum_{n=1}^\infty e^{-nX_r}\nonumber\\
&=&(-)^{q-j-1}e^{X_r}D^{q-j-1} (e^{X_r}-1)^{-1}\qquad (|\arg\,a|<\pi/\lambda),\label{e44}
\end{eqnarray}
with $D\equiv d/dX_r$. The coefficients $c_j$ ($0\leq j\leq \mu-1$) are obtained by recursive solution of (\ref{e43a}) and (\ref{e43}) and are given below\footnote{The values of $c_1$, $c_2$ and $c_3$ follow from (\ref{e42}).} for $\mu=1, 2, \ldots , 5$:
\bee\label{e45}
\begin{array}{ll}
\vspace{0.25cm}

\mu=1: & c_0=1 \\ \vspace{0.25cm}

\mu=2: & c_0=1,\ \ c_1=\gamma-\lambda-1 \\ \vspace{0.25cm}

\mu=3: & c_0=1,\ \ c_1=2\gamma-3\lambda+3,\ \ c_2=(1+\gamma-\lambda)(1+\gamma-2\lambda) \\ 

\mu=4: & c_0=1,\ \ c_1=3\gamma-6\lambda+6,\ \ c_2=3\gamma^2+9\gamma+7-18\lambda-12\lambda\gamma+11\lambda^2,\\\vspace{0.25cm}

       & c_3=(1+\gamma-\lambda)(1+\gamma-2\lambda)(1+\gamma-3\lambda) \\

\mu=5: & c_0=1,\ \ c_1=4\gamma-10\lambda+10,\ \ c_2=6\gamma^2+24\gamma+25-60\lambda-30\lambda\gamma+35\lambda^2,\\
       & c_3=(2\gamma-5\lambda+3)(2\gamma^2+6\gamma+5-15\lambda-10\lambda\gamma+10\lambda^2),\\
       & c_4=(1+\gamma-\lambda)(1+\gamma-2\lambda)(1+\gamma-3\lambda)(1+\gamma-4\lambda). \\
\end{array}
\ee

The generalised Mathieu series when $\mu$ is an integer can therefore be expressed by the following closed-form evaluation.
\begin{theorem}$\!\!\!.$
Let the parameters $\gamma$, $\lambda$ and $\omega_r$ be as in Theorem 1. Let $\mu=q$, where $q$ is a positive integer,
and $\delta=\lambda q-\gamma$, $\vartheta=q+\gamma$. Then with $X_r=2\pi ae^{i\omega_r}$, the generalised Mathieu series has the closed-form evaluation
\bee\label{e46}
S_{q,\gamma}(a;\lambda)=\frac{\g(\frac{\gamma+1}{\lambda}) \g(q-\frac{\gamma+1}{\lambda})}{\lambda \g(q) a^{\delta-1}}+H_{q,\gamma}(a;\lambda)+\frac{(-)^m}{\g(q)} \bl(\frac{\pi}{p}\br)^q a^{q-\delta} \sum_{r=0}^{p-1} E_r(a),
\ee
where $H_{q,\gamma}(a;\lambda)$ is defined in (\ref{e25}) and 
\bee\label{e47}
E_r(a)=e^{-X_r+i\vartheta\omega_r} \sum_{j=0}^{q-1} (-)^jc_j X_r^{-j} K_j(X_r;q)
\ee
with the sums $K_j(X_r;q)$ expressed as derivatives of the exponential term in (\ref{e44}). The coefficients $c_j$
$(0\leq j\leq q-1)$ are given in (\ref{e45}) for $q\leq 5$.
\end{theorem}

\vspace{0.6cm}

\begin{center}
{\bf 5.\ Numerical results and concluding remarks}
\end{center}
\setcounter{section}{5}
\setcounter{equation}{0}
\renewcommand{\theequation}{\arabic{section}.\arabic{equation}}
We present some examples of the large-$a$ expansion of $S_{\mu,\gamma}(a;\lambda)$ given in Sections 3 and 4 to demonstrate numerically the validity of our results. 
\vspace{0.2cm}

\noindent {\it Example 1.}\ \ 
We select $\lambda=4$ ($p=2$, $N=1$), so that for $\gamma=2m$ and $\mu>0$ we obtain from Theorem 2
\bee\label{e51}
S_{\mu,\gamma}(a;4)-\frac{\g(\frac{\gamma+1}{4})
\g(\mu-\frac{\gamma+1}{4})}{4\g(\mu)\,a^{4\mu-\gamma-1}}-H_{\mu,\gamma}(a;4)=\frac{(-)^m (\pi/2)^\mu}{\g(\mu) \,a^{3\mu-\gamma}}\,E_0^*(a),
\ee
as $a\ra+\infty$, where
\[E_0^*(a)=2e^{-X/\surd 2}\bl\{\sum_{j=0}^{M-1} (-)^j c_j X^{-j} K_j^*(X;\f{1}{4}\pi)+O(X^{-M})\br\},\] 
and 
\[K_j^*(X;\f{1}{4}\pi)=\sum_{n=1}^\infty \frac{e^{-(n-1)X/\surd 2}}{n^{1-\mu+j}}\,\cos\,\bl[\frac{nX}{\surd 2}+(j-\mu-\gamma)\frac{\pi}{4}\br]\]
with $X=2\pi a$ and the coefficients $c_j\equiv c_j(\mu,\gamma)$. The leading term on the right-hand side of (\ref{e51}) is easily seen to be given by
\[\frac{2(-)^m (\pi/2)^\mu}{\g(\mu)\,a^{3\mu-\gamma}}\,e^{-X/\surd 2}\,\cos\,\bl[\frac{X}{\surd 2}-(\mu+\gamma)\frac{\pi}{4}\br]\qquad (a\ra+\infty).\]
This last approximation agrees with that obtained in \cite{T} using different methods.

For numerical comparison, we take $\mu=\f{5}{4}$ and three values of $\gamma=0, \pm 2$. Then, we have from (\ref{e51}) the expansion 
\bee\label{e52}
{\hat S}:=S_{\frac{5}{4},\gamma}(a;4)-\frac{\g(\frac{\gamma+1}{4}) \g(1-\frac{\gamma}{4})}{\g(\f{1}{4})\,a^{4-\gamma}}-H_{\frac{5}{4},\gamma}(a;4)\hspace{6cm}\]
\[\sim 
\frac{(-)^m2^{7/4} \pi^{5/4}}{\g(\f{1}{4})\,a^{15/4-\gamma}}\,e^{-X/\surd 2}
\sum_{j=0}^\infty \frac{(-)^j c_j}{X^j} \sum_{n=1}^\infty \frac{e^{-(n-1)X/\surd 2}}{n^{j-1/4}}\,\cos\,\bl[\frac{nX}{\surd 2}+(j-\frac{5}{4}-\gamma)\frac{\pi}{4}\br],
\ee
as $a\ra+\infty$, where from (\ref{e25}) the algebraic expansion is
\[H_{\frac{5}{4},\gamma}(a;4)=0\ \ (\gamma=2), \quad -\frac{1}{2a^5}\ \ (\gamma=0),\quad \frac{\pi^2}{6a^5}\ \ (\gamma=-2).\]
The coefficients $c_j$ (with $c_0=1$) for the above three values of $\gamma$ are tabulated in Table 1 for $1\leq j\leq 10$.
In Table 2, we show the absolute relative error in the computation of ${\hat S}$ (defined as the left-hand side of (\ref{e52})) for different $\gamma$ and truncation index $j$ using the expansion (\ref{e52}) when $a=5$. The corresponding value of ${\hat S}$ is indicated at the head of each column. 
In these computations the sum over $n$ was evaluated to an accuracy commensurate with the overall level of precision.
The index $j=12$ corresponds approximately to optimal truncation of $E_0^*(a)$; that is, truncation at or near the least term in absolute value.
\begin{table}[t]
\caption{\footnotesize{The absolute relative error in the computation of ${\hat S}$ from (\ref{e52}) for different $\gamma$ and truncation index $j$ when $\mu=5/4$, $\lambda=4$ and $a=5$.}}
\begin{center}
\begin{tabular}{|r|c|c|c|}
\hline
&&&\\[-0.3cm]
\mcol{1}{|c|}{} & \mcol{1}{c|}{$\gamma=0$} &\mcol{1}{c|}{$\gamma=2$}& \mcol{1}{c|}{$\gamma=-2$} \\
\mcol{1}{|c|}{$j$} & \mcol{1}{c|}{${\hat S}=-1.54766(-12)$} &\mcol{1}{c|}{${\hat S}=-3.59325(-11)$}& \mcol{1}{c|}{${\hat S}=+5.75174(-14)$} \\
[.1cm]\hline
&&&\\[-0.25cm]
0 & $1.980(-02)$ & $1.178(-04)$ & $3.237(-04)$ \\
1 & $3.378(-04)$ & $1.568(-04)$ & $1.530(-03)$ \\
2 & $1.102(-07)$ & $1.086(-05)$ & $2.066(-04)$ \\
3 & $3.649(-07)$ & $1.758(-07)$ & $1.920(-05)$ \\
4 & $2.697(-08)$ & $4.874(-09)$ & $3.210(-07)$ \\
6 & $2.785(-11)$ & $1.265(-09)$ & $1.437(-07)$ \\
8 & $1.220(-11)$ & $3.334(-12)$ & $1.363(-09)$ \\
10& $9.421(-14)$ & $1.497(-12)$ & $9.459(-10)$ \\
12& $2.607(-14)$ & $1.277(-14)$ & $2.768(-11)$ \\
[.2cm]\hline
\end{tabular}
\end{center}
\end{table}
\vspace{0.2cm}

\noindent {\it Example 2.}\ \ 
In Theorem 3, we first consider the case $\mu=1$ where, from (\ref{e47}), 
\[e^{-i\vartheta\omega_r}E_r(a)=e^{-X_r}\, K_0(X_r;1)=\frac{e^{-X_r}}{1-e^{-X_r}}=\frac{1}{2}\coth (\pi ae^{i\omega_r})-\frac{1}{2}.\]
Some straightforward algebra shows, when $\vartheta=1+2m$, that
\[\sum_{r=0}^{p-1}e^{i\vartheta\omega_r}=e^{\pi i\vartheta(p-1)/(2p)} \sum_{r=0}^{p-1} e^{-\pi i\vartheta r/p}=\frac{\sin\,(\pi\vartheta/2)}{\sin\,(\pi\vartheta/\lambda)}
=\frac{(-)^m}{\sin\,(\pi\vartheta/\lambda)}.\] 
It then follows from (\ref{e46}) (where the first term in (\ref{e46}) involving $a^{1-\delta}$ cancels with the contribution from the above finite sum) that
\[S_{1,2m}(a;2p)=\sum_{n=1}^\infty \frac{n^{2m}}{n^{2p}+a^{2p}}\hspace{7cm}\]
\bee\label{e53}
\hspace{1.2cm}=H_{1,2m}(a;2p)+\frac{(-)^m \pi}{2p a^{2p-2m-1}}\, \sum_{r=0}^{p-1} e^{i\vartheta\omega_r}\coth (\pi ae^{i\omega_r})~.
\ee
Series of this type have been expressed as infinite sums of Riemann zeta functions in \cite{CHP}.

In the case $\gamma=0$, $\lambda=2$, we have $m=0$, $p=1$, $\omega_0=0$ and $H_{1,0}(a;2)=-1/(2a^2)$. Then (\ref{e53}) yields the well-known result
\[S_{1,0}(a;2)=\sum_{n=1}^\infty \frac{1}{n^2+a^2}=\frac{\pi}{2a} \coth \pi a-\frac{1}{2a^2}~.\]
When $\mu=3$, $\gamma=\lambda=2$, we have $\omega_0=0$ and from (\ref{e47})
\[E_0(a)=\sum_{j=0}^2 c_jX^{-j}\,D^{2-j} (e^X-1)^{-1},\qquad X=2\pi a.\]
The coefficients for these values of $\mu$, $\gamma$ and $\lambda$ are, from (\ref{e45}), found to be $c_0=1$, $c_1=1$ and $c_2=-1$, whence
\[S_{3,2}(a;2)=\sum_{n=1}^\infty \frac{n^2}{(n^2+a^2)^3}=\frac{\pi}{16a^3}-\frac{\pi^3}{2a}\,\frac{e^{-X}}{1-e^{-X}}\bl\{\frac{(1+e^{-X})}{(1-e^{-X})^2}+\frac{1}{X(1-e^{-X})}-\frac{1}{X^2}\br\}.\] 
Similarly, if $\gamma=2$, $\lambda=4$ we have $\omega_0=\f{1}{4}$, $\omega_1=-\f{1}{4}$ and $c_0=1$, $c_1=-5$ and $c_2=5$. Then, for $a>0$, we obtain
\[S_{3,2}(a;4)=\sum_{n=1}^\infty \frac{n^2}{(n^4+a^4)^3}
=\frac{5\pi\sqrt{2}}{128a^9}\hspace{6cm}\]
\[\hspace{2cm}+\frac{\pi^3}{8a^7}\,\Re \bl\{
\frac{e^{-X_0+\pi i/4}}{1-e^{-X_0}}\bl[
\frac{(1+e^{-X_0})}{(1-e^{-X_0})^2}+\frac{5}{X_0(1-e^{-X_0})}+\frac{5}{X_0^2}\br]\br\},\quad X_0=2\pi ae^{\pi i/4}.\]
\vspace{0.2cm} 

Finally, we remark that the asymptotic expansion of the alternating version of (\ref{e12}) can be deduced by making use of the identity
\bee\label{e56}
{\tilde S}_{\mu,\gamma}(a;\lambda):=\sum_{n=1}^\infty \frac{(-)^{n-1} n^\gamma}{(n^\lambda+a^\lambda)^\mu}=S_{\mu,\gamma}(a;\lambda)-2^{1-\delta} S_{\mu,\gamma}(\fs a;\lambda).
\ee
Substitution of the expansion for $S_{\mu,\gamma}(a;\lambda)$ in (\ref{e38}) into (\ref{e56}) leads to the introduction of the alternating analogues ${\tilde H}_{\mu,\gamma}(a;\lambda)$ and ${\tilde E}_r(a)$ of the algebraic and exponential expansions given by (with $k^*$ defined in (\ref{e23a}))
\bee\label{e54}
{\tilde H}_{\mu,\gamma}(a;\lambda)=\left\{\begin{array}{ll} \!\!0 & \ \ (m=1, 2, \ldots) \\
\displaystyle{
\!\!\frac{a^{-\lambda\mu}}{\g(\mu)}\sum_{k=0}^{k^*}\frac{(-)^k \g(\mu+k)}{k! \,a^{\lambda k}}\,(1-2^{2m+1+\lambda k})\,\zeta(2m-2kp)} & \ \ (m=0, -1, -2, \ldots) \end{array}\right.
\ee
and  
\bee\label{e55}
{\tilde E}_r(a)=e^{-X_r/2+i\vartheta\omega_r}\bl\{\sum_{j=0}^{M-1}(-)^jc_j X_r^{-j} {\tilde K}_j(X_r;\mu)+O(X_r^{-M})\br\},
\ee
for $0\leq r\leq p-1$, where
\begin{eqnarray*}
{\tilde K}_j(X_r;\mu)&=&-e^{X_r/2}\{e^{-X_r}K_j(X_r;\mu)-2^{1-\mu+j} e^{-X_r/2} K_j(\fs X_r;\mu)\}\\
&=&-e^{X_r/2} \sum_{n=1}^\infty \bl\{\frac{e^{-nX_r}}{n^{1-\mu+j}}-\frac{e^{-nX_r/2}}{(\fs n)^{1-\mu+j}}\br\}\\
&=&\sum_{n=1}^\infty 
\frac{e^{-(n-1)X_r}}{(n\!-\!\fs)^{1-\mu+j}}~.
\end{eqnarray*}
We then have the following theorem.
\begin{theorem}$\!\!\!.$
Let the parameters $\mu$, $\gamma$, $\lambda$ and the quantities $\delta$, $\vartheta$,  $\omega_r$ be as in Theorem 1. Then the expansion of the alternating series is
\[{\tilde S}_{\mu,\gamma}(a;\lambda)={\tilde H}_{\mu,\gamma}(a;\lambda)+\frac{(-)^{m-1}}{\g(\mu)} 
\bl(\frac{\pi}{p}\br)^\mu a^{\mu-\delta} \sum_{r=0}^{p-1} {\tilde E}_r(a)\]
as $|a|\ra\infty$ in $|\arg\,a|<\pi/\lambda$, where ${\tilde H}_{\mu,\gamma}(a;\lambda)$ and ${\tilde E}_r(a)$ are defined in (\ref{e54}) and (\ref{e55}).
\end{theorem}


\vspace{0.6cm}

\begin{thebibliography}{99}
\bibitem{Br}
B.L.J. Braaksma, Asymptotic expansions and analytic continuations for a class of Barnes integrals, Compos. Math. {\bf 15} (1963) 239--341.

\bibitem{M}
E.L. Mathieu, {\it Trait\'e de Physique Math\'ematique. VI--VII: Th\'eorie de l'Elasticit\'e des Corps Solides} (Part 2), Gauthier-Villars, Paris, 1890.

\bibitem{DLMF}
F.W.J. Olver, D.W. Lozier, R.F. Boisvert and C.W. Clark (eds.),    
{\it NIST Handbook of Mathematical Functions}, Cambridge University Press, Cambridge, 2010.

\bibitem{P}
R.B. Paris, The asymptotic expansion of a generalised Mathieu series, Appl. Math. Sci. {\bf 125} (2013) 6209--6216.

\bibitem{P2}
R.B. Paris, The asymptotic expansion of a generalisation of the Euler-Jacobi series, Eur. J. Pure Appl. Math. [in press]  arXiv:1503.07329 (2015).

\bibitem{PK} 
R.B. Paris and D. Kaminski,  {\it Asymptotics and Mellin-Barnes Integrals}, 
Cambridge University Press, Cambridge, 2001.

\bibitem{CHP}
C.H. Picard, On some series formed by values of the Riemann zeta function. arXiv: 1511.04720 (2015).

\bibitem{PST}
T. Pog\'any, H.M. Srivastava and Z. Tomovski, Some families of Mathieu {\bf a}-series and alternating Mathieu {\bf a}-series Appl. Math. Comp. {\bf 173} (2006) 69--108.

\bibitem{ECT}
E.C. Titchmarsh, {\it The Theory of the Riemann Zeta-Function}, revised by D.R. Heath-Brown, Oxford University Press, Oxford, 1986.

\bibitem{T}
K. Tsouvalas, private communication.


\bibitem{Z}
V.P. Zastavnyi, Asymptotic expansions of several series and their application, Ukrainian Math. Bull. {\bf 6} (2009) 549--569.
\end{thebibliography}
\end{document}